\documentclass{article}
\let\OLDthebibliography\thebibliography
\renewcommand\thebibliography[1]{
  \OLDthebibliography{#1}
  \setlength{\parskip}{0pt}
  \setlength{\itemsep}{0pt plus 0.3ex}}
  
\title{Principal Branches of Inverse Trigonometric and Inverse Hyperbolic Functions}
\author{Kevin M. Dempsey}
\date{\today}

\newcommand{\D}{\displaystyle}
\newcommand{\rmd}{\mathrm{d}}
\newcommand{\rmi}{\mathrm{i}}
\newcommand{\XDom}{\mathbb{C}\setminus(-1,1)}
\newcommand{\YDom}{\mathbb{C}\setminus(-\infty,-1)\cup(1,\infty)}
\newcommand{\XoneDom}{\mathbb{C}\setminus(-\infty,1)}
\newcommand{\XtwoDom}{\mathbb{C}\setminus(-\infty,-1)}
\newcommand{\XRange}{\mathbb{C}\setminus(-\rmi,\rmi)}
\newcommand{\arcsincosDom}{\mathbb{C}\setminus(-\infty,-1)\cup(1,\infty)}
\newcommand{\arccscsecDom}{\mathbb{C}\setminus(-1,1)}
\newcommand{\arctanDom}{\mathbb{C}\setminus(-\rmi\infty,-\rmi]\cup[\rmi,\rmi\infty)}
\newcommand{\arccotDom}{\mathbb{C}\setminus(-\rmi,\rmi)}
\newcommand{\arcsinhDom}{\mathbb{C}\setminus(-\rmi\infty,-\rmi)\cup(\rmi,\rmi\infty)}
\newcommand{\arccoshDom}{\mathbb{C}\setminus(-\infty,1)}
\newcommand{\arccschDom}{\mathbb{C}\setminus(-\rmi,\rmi)}
\newcommand{\arcsechDom}{\mathbb{C}\setminus(-\infty,0]\cup(1,\infty)}
\newcommand{\arctanhDom}{\mathbb{C}\setminus(-\infty,-1]\cup[1,\infty)}
\newcommand{\arccothDom}{\mathbb{C}\setminus[-1,1]}

\usepackage{amssymb}
\usepackage{amsmath}
\DeclareMathOperator{\re}{Re}
\DeclareMathOperator{\sign}{sign}
\DeclareMathOperator{\arccsc}{arccsc}
\DeclareMathOperator{\arcsec}{arcsec}
\DeclareMathOperator{\arccot}{arccot}
\DeclareMathOperator{\arcsinh}{arcsinh}
\DeclareMathOperator{\arccosh}{arccosh}
\DeclareMathOperator{\arccsch}{arccsch}
\DeclareMathOperator{\arcsech}{arcsech}
\DeclareMathOperator{\arctanh}{arctanh}
\DeclareMathOperator{\arccoth}{arccoth}
\DeclareMathOperator{\X}{X}
\DeclareMathOperator{\Y}{Y}
\DeclareMathOperator{\Z}{Z}

\begin{document}
\maketitle

\begin{abstract}
We develop principal branches for three key square root functions and for the inverse trigonometric and inverse hyperbolic functions.
The three square root branches are integral to defining the inverse function branches, their derivatives, and their antiderivatives.
Complex analysis is used to turn the definitions of the principal branches into concrete expressions.
We take the standard reference in this area to be the NIST Digital Library of Mathematical Functions (DLMF).
We adopt the notation for, and the definitions of, the principal branches of the inverse functions in the DLMF.
Our goal is to widen the scope of the results in the DLMF while at the same time lowering the complex variables burden on the average DLMF user.
\end{abstract}

\section{Introduction}

With the advent of large scale computing machines over seventy years ago, the need for mathematical tables was questioned. The publication of such tables dates at least back to 1614 and the tables of natural logarithms by John Napier\cite{Napier}. After some consideration, the National Science Foundation asked the National Bureau of Standards (NBS) to modernize the existing classical tables in a way that would complement the changing times. In 1964, the NBS published the result of their efforts, a handbook of mathematical functions edited by Milton Abramowitz and Irene Stegun\cite{A&S, Grier}. The handbook came to be known as Abramowitz and Stegun (A\&S). Any fears that such a work was past its time were quelled. The handbook had a wider circulation and was cited more than any other publication in the first 100 years of the institute. The NBS was founded in 1901 and was renamed the National Institute of Standards and Technology (NIST) in 1988. North of a thousand pages, A\&S largely succeeded in defining, and standardizing the notation for, special functions. Thirty plus years after its publication, after numerous printings with corrections and published errata, the decision was made at NIST to bring A\&S into the 21st century.
In 2010, after ten years in the making, NIST made available on a publicly available website, the Digital Library of Mathematical Functions (DLMF)\cite{DLMF}. The DLMF is an online counterpart to the printed A\&S. Moreover, it is continually updated and validated by NIST. No longer is there a need to publish errata. The DLMF website documents the versions and lists the updates. By design, the DLMF is easy to navigate. It is a free handbook that you do not have to carry with you.

The branch cuts for the inverse trigonometric functions and for the inverse hyperbolic functions are depicted schematically in the DLMF\cite[Fig.~4.23.1,~4.37.1]{DLMF}.
As far as the branch cut intervals being open or closed at the endpoints, we adhere to the conventions adopted by Wolfram\cite{Weisstein}.

We caution the reader that the choice of principal branch for the inverse cotangent is not universal. The principal branch of the inverse cotangent in the DLMF, and here, is the branch returned by Mathematica\cite{Mathematica}, but it is not the branch returned by Maple\cite{Maple}, for example.

Throughout this article, $z=x+\rmi y$ is a complex number with real part $x$ and imaginary part $y$, $\ln{z}$ is the principal branch of the natural logarithm\cite[Eq.~4.2.2]{DLMF}, and $\sign(x)$ is the sign function\cite{DLMF}
$$%\begin{equation}
\sign(x)=\left\{\!\!\begin{array}{rc} -1, & x<0, \\ 0, & x=0, \\ 1, & x>0. \end{array} \right.
$$%\end{equation}
We note that
$$%\begin{equation}
\sign(x)\arccosh(\sqrt{1+x^2}) = \arcsinh{x} \quad (x\in\mathbb{R}).
$$%\end{equation}

\section{Principal branch $\X_1(z)$ of $\sqrt{z-1}$ on $\XoneDom$}
The square root function $\sqrt{z-1}$ has branch points at $z=1$ and $\infty$.
We define the principal branch $\X_1(z)$ of $\sqrt{z-1}$, the analytic continuation of $\sqrt{x-1}$ from $[1,\infty)$ to the cut plane $\XoneDom$,
with $r_1=|z-1|$, by\cite[p.~260]{Johnson}
$$%\begin{equation}
\X_1(z) = \sqrt{\frac{r_1+(x-1)}{2}}+\rmi\sign(y)\sqrt{\frac{r_1-(x-1)}{2}}.
$$%\end{equation}
Clearly, $|\X_1(z)|=\sqrt{r_1}$.
On the $x$-axis away from the cut, we confirm that
$$%\begin{equation}
\X_1(x) = \sqrt{x-1} \quad (x\ge1).
$$%\end{equation}
On the $x$-axis on the cut, $\X_1(z)$ is two-valued,
$$%\begin{equation}
\X_1(x\pm\rmi0) = +0\pm\rmi\sqrt{1-x} \quad (x<1).
$$%\end{equation}
On the $y$-axis, with $r=\sqrt{1+y^2}$,
$$%\begin{equation}
\X_1(\rmi y) = \sqrt{\frac{r-1}{2}}+\rmi\sign(y)\sqrt{\frac{r+1}{2}} \quad (y^2>0).
$$%\end{equation}
We see that $\X_1(z)$ maps the cut plane $\XoneDom$ to the closed right half plane.
The upper half plane maps to the first quadrant, and the lower half plane maps to the fourth quadrant.We verify that $\X_1^2(z)=z-1$ and $\X'_1(z)=\frac{1}{2\X_1(z)}$.

\section{Principal branch $\X_2(z)$ of $\sqrt{z+1}$ on $\XtwoDom$}
The square root function $\sqrt{z+1}$ has branch points at $z=-1$ and $\infty$.
We define the principal branch $\X_2(z)$ of $\sqrt{z+1}$, the analytic continuation of $\sqrt{x+1}$ from $[-1,\infty)$ to the cut plane $\XtwoDom$,
with $r_2=|z+1|$, by\cite[p.~260]{Johnson}
$$%\begin{equation}
\X_2(z) = \sqrt{\frac{r_2+(x+1)}{2}}+\rmi\sign(y)\sqrt{\frac{r_2-(x+1)}{2}}.
$$%\end{equation}
Clearly, $|\X_2(z)|=\sqrt{r_2}$.
On the $x$-axis away from the cut, we confirm that
$$%\begin{equation}
\X_2(x) = \sqrt{x+1} \quad (x\ge-1).
$$%\end{equation}
On the $x$-axis on the cut, $\X_2(z)$ is two-valued,
$$%\begin{equation}
\X_2(x\pm\rmi0) = +0\pm\rmi\sqrt{-1-x} \quad (x<-1).
$$%\end{equation}
On the $y$-axis, again with $r=\sqrt{1+y^2}$,
$$%\begin{equation}
\X_2(\rmi y) = \sqrt{\frac{r+1}{2}}+\rmi\sign(y)\sqrt{\frac{r-1}{2}} \quad (y^2>0).
$$%\end{equation}
We see that $\X_2(z)$ maps the cut plane $\XtwoDom$ to the closed right half plane.
The upper half plane maps to the first quadrant, and the lower half plane maps to the fourth quadrant.
We verify that $\X_2^2(z)=z+1$ and $\X'_2(z)=\frac{1}{2\X_2(z)}$.

We note that for $y^2>0$,
$$%\begin{equation}
\X_1(-z)=-\rmi\sign(y)\X_2(z), \quad \X_2(-z)=-\rmi\sign(y)\X_1(z).
$$%\end{equation}

\section{Principal branch $\X(z)$ of $\sqrt{z^2-1}$ on $\XDom$}
The square root function $\sqrt{z^2-1}$ has branch points at $\pm1$.
On the $x$-axis away from the cuts, we see that
$$%\begin{equation}
\X_1(x)\X_2(x) = \sqrt{x^2-1} \quad (x\ge1).
$$%\end{equation}
By expanding and simplifying $\X_1(z)\X_2(z)$ for $y\ne 0$, we discover that the principal branch of $\sqrt{z^2-1}$ on the cut plane $\XDom$ is given by
$$%\begin{equation}
\X(z) = \beta\sqrt{\alpha^2-1}+\rmi\sign(y)\alpha\sqrt{1-\beta^2},
$$%\end{equation}
where $\alpha=(r_1+r_2)/2$ and $\beta=(r_2-r_1)/2$ with $\alpha\ge 1$, $\beta\in[-1,1]$ ($\beta>0$ for $x>0$ and $\beta<0$ for $x<0$).
We note that $|\X(z)|=\sqrt{r_1r_2}$.
Recall that $r_1=|z-1|=\sqrt{(x-1)^2+y^2}$ and $r_2=|z+1|=\sqrt{(x+1)^2+y^2}$.
Since $r_1(-z)=r_2(z)$ and $r_2(-z)=r_1(z)$, $\alpha(-z)=\alpha(z)$, $\beta(-z)=-\beta(z)$ and $\X(-z)=-\X(z)$.
On the $x$-axis away from the cut, $\alpha=|x|$, $\beta=\sign(x)$, and
$$%\begin{equation}
\X(x)=\sign(x)\sqrt{x^2-1} \quad (x^2\ge1).
$$%\end{equation}
On the $x$-axis, on the cut, $\alpha=1$, $\beta=x$, and $\X(z)$ is two-valued,
$$%\begin{equation}
\X(x\pm\rmi 0)=\sign(x)0\pm\rmi\sqrt{1-x^2} \quad (x^2<1).
$$%\end{equation}
We note that $\X(0\pm\rmi 0)=\pm\rmi$.
On the $y$-axis, $\alpha=\sqrt{1+y^2}$, $\beta=0$, and $\X(z)$ is purely imaginary,
$$%\begin{equation}
\X(\rmi y)=\rmi\sign(y)\sqrt{1+y^2} \quad (y^2>0).
$$%\end{equation}
We see that $\X(z)$ maps the cut plane $\XDom$ to the cut plane $\XRange$.
Away from the axes, $z$ and $\X(z)$ are in the same quadrant.
We verify that $\X(z)^2=z^2-1$ and $\X'(z)=z/\X(z)$, and, for completeness,
we note that the partial derivatives of $\alpha$ and $\beta$ with respect to $x$ are
$$%\begin{equation}
\alpha_x=\frac{\beta(\alpha^2-1)}{\alpha^2-\beta^2},\quad \beta_x=\frac{\alpha(1-\beta^2)}{\alpha^2-\beta^2}.
$$%\end{equation}

\section{Principal branch $\Y(z)$ of $\sqrt{1-z^2}$ on $\YDom$}
If we look back at the principal branch $\X_2(z)$ of $\sqrt{z+1}$, we note that
$$%\begin{equation}
\X_2(-x)\X_2(x)=\sqrt{1-x^2} \quad (x^2\le 1).
$$%\end{equation}
For $y\ne 0$, $\X_2(-z)\X_2(z)=-\rmi\sign(y)\X_1(z)\X_2(z)$ and we find that the principal branch of $\sqrt{1-z^2}$ on the cut plane $\YDom$ is given by
$$%\begin{equation}
\Y(z) = \alpha\sqrt{1-\beta^2}-\rmi\sign(y)\beta\sqrt{\alpha^2-1}.
$$%\end{equation}
In this case, $|\Y(z)|=\sqrt{r_1r_2}$ and $\Y(-z)=\Y(z)$.
On the $x$-axis away from the cut, $\alpha=1$, $\beta=x$, and we confirm that
$$%\begin{equation}
\Y(x) = \sqrt{1-x^2} \quad (x^2\le 1).
$$%\end{equation}
On the $x$-axis on the cut, $\alpha=|x|$, $\beta=\sign(x)$, and $\Y(z)$ is two-valued
$$%\begin{equation}
\Y(x\pm\rmi 0) = +0\mp\rmi\sign(x)\sqrt{x^2-1} \quad (x^2>1).
$$%\end{equation}
On the $y$-axis, $\alpha=\sqrt{1+y^2}$, $\beta=0$, and $\Y(z)$ is real,
$$%\begin{equation}
\Y(\rmi y) = \sqrt{1+y^2} \quad (y\in\mathbb{R}).
$$%\end{equation}
We see that $\Y(z)$ maps the cut plane $\YDom$ to the closed right half plane.
Away from the axes, the first and third quadrants map under $\Y(z)$ to the fourth quadrant while the second and fourth quadrants map to the first quadrant.
We verify that $\Y(z)^2=1-z^2$ and $\Y'(z)=-z/\Y(z)$.

\section{Principal branches of the inverse sine and inverse cosine}
The inverse sine and inverse cosine functions have branch points at $\pm1$ and $\infty$\,\cite{Wolfram}.
The principal branches $\arcsin{z}$ and $\arccos{z}$ are defined on the cut plane $\arcsincosDom$ by\cite[\S4.23.1--2]{DLMF}
$$%\begin{equation}
\begin{aligned}
\arcsin{z} &= \int_{0}^{z}\frac{\rmd t}{\sqrt{1-t^2}} = \int_{0}^{z}\frac{\rmd t}{\Y(t)}, \\
\arccos{z} &= \int_{z}^{1}\frac{\rmd t}{\sqrt{1-t^2}} = \int_{z}^{1}\frac{\rmd t}{\Y(t)},
\end{aligned}
$$%\end{equation}
where the principal branch $\Y(z)$ defined above is the analytic continuation of $\sqrt{1-x^2}$ from $[-1,1]$ to the cut plane $\arcsincosDom$.
If we let $s(z)=\arcsin{\beta}+\rmi\sign(y)\arccosh{\alpha}$ we find that $s'(z)=1/\Y(z)$.
At $z=0$, $\alpha=1$ and $\beta=0$, while at $z=1$, $\alpha=1$ and $\beta=1$.
It follows that
$$%\begin{equation}
\begin{aligned}
\arcsin{z} &= s(z)-s(0) = \arcsin{\beta}+\rmi\sign(y)\arccosh{\alpha}, \\
\arccos{z} &= s(1)-s(z) = \arccos{\beta}-\rmi\sign(y)\arccosh{\alpha}.
\end{aligned}
$$%\end{equation}
These expressions are in Johnson\cite[p.~264]{Johnson}.
The DLMF\cite[Eq.~4.23.34--35,37--38]{DLMF} expresses these equations in terms of the natural logarithm.
The second equation is $\arccos{z}=\pi/2-\arcsin{z}$.
We see that $\arcsin(-z)=-\arcsin{z}$ and $\arccos(-z)=\pi-\arccos{z}$.
These reflection formulae are in the DLMF\cite[Eq.~4.23.10--11,16]{DLMF}.
On the $x$-axis away from the cut, $z=x\in[-1,1]$, $\alpha=1$, $\beta=x$, and we confirm that $\arcsin{z}=\arcsin{x}$ and $\arccos{z}=\arccos{x}$.
The principal branches $\arcsin{z}$ and $\arccos{z}$ are the analytic continuations of $\arcsin{x}$ and $\arccos{x}$ from $[-1,1]$ to the cut plane $\arcsincosDom$.
On the $x$-axis on the cut, $z=x\pm\rmi 0 (x^2>1)$, $\alpha=|x|$, $\beta=\sign(x)$, and $\arcsin{z}$ and $\arccos{z}$ are two-valued,
$$%\begin{equation}
\begin{aligned}
\arcsin(x\pm\rmi 0) &= \sign(x)\frac{\pi}{2}\pm\rmi\arccosh{|x|} & (x^2>1), \\
\arccos(x\pm\rmi 0) &= \arccos[\sign(x)]\mp\rmi\arccosh{|x|} & (x^2>1).
\end{aligned}
$$%\end{equation}
The DLMF\cite[Eq.~4.23.20--21,24--25]{DLMF} expresses these results in terms of the natural logarithm.
On the $y$-axis, $\alpha=\sqrt{1+y^2}$, $\beta=0$ and
$$%\begin{equation}
\begin{aligned}
\arcsin(\rmi y) &= \rmi\arcsinh{y}               & (y\in\mathbb{R}), \\
\arccos(\rmi y) &= \frac{\pi}{2}-\rmi\arcsinh{y} & (y\in\mathbb{R}).
\end{aligned}
$$%\end{equation}
We see that the cut plane $\arcsincosDom$ is mapped to $[-\pi/2,\pi/2]\times\mathbb{R}$ by $\arcsin{z}$ and to $[0,\pi]\times\mathbb{R}$ by $\arccos{z}$.
From the development above
$$%\begin{equation}
\begin{aligned}
\frac{\rmd}{\rmd z}\arcsin{z} &=  \frac{1}{\Y(z)}, & \int\arcsin{z}\,\rmd z &= z\arcsin{z}+\Y(z)+C, \\
\frac{\rmd}{\rmd z}\arccos{z} &= -\frac{1}{\Y(z)}, & \int\arccos{z}\,\rmd z &= z\arccos{z}-\Y(z)+C.
\end{aligned}
$$%\end{equation}

\section{Principal branches of the inverse cosecant and inverse secant}
The inverse cosecant and inverse secant functions have branch points at $0$ and $\pm1$\,\cite{Wolfram}.
The principal branches $\arccsc{z}$ and $\arcsec{z}$ are defined on the cut plane $\arccscsecDom$ by\cite[Eq.~4.23.7--8]{DLMF}
$$%\begin{equation}
\begin{aligned}
\arccsc{z} &= \arcsin(1/z), \\
\arcsec{z} &= \arccos(1/z).
\end{aligned}
$$%\end{equation}
It follows that $\arccsc(-z)=-\arccsc{z}$, $\arcsec(-z)=\pi-\arcsec{z}$ and $\arcsec{z}=\pi/2-\arccsc{z}$.
These reflection formulae are in the DLMF\cite[Eq.~4.23.13--14,17]{DLMF}.
On the $x$-axis away from the cut, $x^2\ge1$ and we confirm that $\arccsc{z}=\arccsc{x}$ and $\arcsec{z}=\arcsec{x}$.
The principal branches $\arccsc{z}$ and $\arcsec{z}$ are the analytic continuations of $\arccsc{x}$ and $\arcsec{x}$ from $(-\infty,-1]\cup[1,\infty)$ to the cut plane $\arccscsecDom$.
On the $x$-axis on the cut away from the branch point at the origin, $0<x^2<1$ and $\arccsc{z}$ and $\arcsec{z}$ are two-valued,
$$%\begin{equation}
\begin{aligned}
\arccsc(x\pm\rmi 0) &= \sign(x)\frac{\pi}{2}\mp\rmi\arcsech|x| & (0<x^2<1), \\
\arcsec(x\pm\rmi 0) &= \arccos[\sign(x)]\pm\rmi\arcsech|x| & (0<x^2<1).
\end{aligned}
$$%\end{equation}
On the $y$-axis,
$$%\begin{equation}
\begin{aligned}
\arccsc(\rmi y) &= -\rmi\arccsch{y}              & (y^2>0), \\
\arcsec(\rmi y) &= \frac{\pi}{2}+\rmi\arccsch{y} & (y^2>0).
\end{aligned}
$$%\end{equation}
We see that the cut plane $\arccscsecDom$ is mapped to $[-\pi/2,\pi/2]\times\mathbb{R}$ by $\arccsc{z}$ and to $[0,\pi]\times\mathbb{R}$ by $\arcsec{z}$.
We verify that
$$%\begin{equation}
\begin{aligned}
\frac{\rmd}{\rmd z}\arccsc{z} &=-\frac{1}{z^2\Y(1/z)}, & \int\arccsc{z}\,\rmd z &= z\arccsc{z}+\arctanh(\Y(1/z))+C, \\
\frac{\rmd}{\rmd z}\arcsec{z} &=\frac{1}{z^2\Y(1/z)},  & \int\arcsec{z}\,\rmd z &= z\arcsec{z}-\arctanh(\Y(1/z))+C.
\end{aligned}
$$%\end{equation}

\section{Principal branches of the inverse tangent and inverse cotangent}
The inverse tangent and inverse cotangent functions have branch points at $\pm\rmi$\,\cite{Wolfram}.
The principal branch $\arctan{z}$ of the inverse tangent is defined on the cut plane $\arctanDom$ by\cite[Eq.~4.23.3]{DLMF}
$$%\begin{equation}
\arctan{z}=\int_0^z\frac{\rmd t}{1+t^2}.
$$%\end{equation}
The principal branch $\arccot{z}$  of the inverse cotangent is defined on the cut plane $\arccotDom$ by\cite[Eq.~4.23.9]{DLMF}
$$%\begin{equation}
\arccot{z}=\arctan(1/z)=\int_{z}^{\infty}\frac{\rmd t}{1+t^2} \quad (x\ne 0).
$$%\end{equation}
On the $x$-axis, we note that
$$%\begin{equation}
\begin{aligned}
\arccot{x} &= \int_{x}^{-\infty}\frac{\rmd t}{1+t^2} & (x<0), \\
\arctan{x}+\arccot{x} &= \sign(x)\frac{\pi}{2} & (x\ne0).
\end{aligned}
$$%\end{equation}
The latter result genaralizes.
Away from the axes, we apply Cauchy's theorem to the sectorial contour that has one arm along the positive or negative $x$-axis, the other arm passing through $z$, and the circular arc in the same quadrant as $z$, to show that
$$%\begin{equation}
\arctan{z}+\arccot{z}=\sign(x)\frac{\pi}{2} \quad (x\ne 0).
$$%\end{equation}
This reflection formula is in the DLMF\cite[Eq.~4.23.18]{DLMF}.
When $x>0$, the rectilinear path of integration,
$$%\begin{equation}
\arccot{z} = \int_{y}^{0}\frac{\rmi\rmd\eta}{1+(x+\rmi\eta)^2}+\int_{x}^{\infty}\frac{\rmd t}{1+t^2},
$$%\end{equation}
leads to a closed-form expression for $\arccot{z}$, namely,
$$%\begin{equation}
\arccot{z} = \frac{1}{2}\left(\arctan\frac{y+1}{x}-\arctan\frac{y-1}{x}\right)+\frac{\rmi}{4}\ln\frac{x^2+(y-1)^2}{x^2+(y+1)^2}\quad(x\ne0).
$$%\end{equation}
When $x<0$, the upper limit of the second integral is $-\infty$ and the result is the same.
The corresponding expression for $\arctan{z}$ for $x\ne0$ follows from
$\arctan{z}=\sign(x)\pi/2-\arccot{z}$.
The expression for $\arctan{z}$ in the DLMF\cite[Eq~4.23.36]{DLMF} holds for $|z|<1$.
We observe that $\arccot(-z)=-\arccot{z}$ and $\arctan(-z)=-\arctan{z}$.
These reflection formulae are in the DLMF\cite[Eq.~4.23.12,15]{DLMF}.
With the aid of 1.625 9 in Gradshteyn \& Ryzhik~\cite{G&R}, and with the aid of a simple sketch in the case when $z=x+\rmi y$ is on the unit circle, we are able to capture the behavior of the real part of $\arccot{z}$ as we go successively from being inside, on, and outside the unit circle ($x\ne 0$ in each case),
$$%\begin{equation}
\re(\arccot{z})
=
\begin{cases}
\D\sign(x)\frac{\pi}{2}-\frac{1}{2}\arctan\frac{2x}{1-x^2-y^2} & (|z|<1), \\
\D\sign(x)\frac{\pi}{4}                                     & (|z|=1), \\
\D\frac{1}{2}\arctan\frac{2x}{x^2+y^2-1}                   & (|z|>1).
\end{cases}
$$%\end{equation}
We see this behavior on the $x$-axis.
As $x\to\pm 0$, $\arccot(x)\to\pm\pi/2$, $\arccot(\pm 1)=\pm\pi/4$, and as $x\to\pm\infty$, $\arccot{x}\to 0$.
On the unit circle away from the imaginary axis,
$$%\begin{equation}
\arccot{z} = \sign(x)\frac{\pi}{4}-\frac{\rmi}{2}\arctanh{y} \quad (y^2<1).
$$%\end{equation}

The principal branch $\arctan{z}$ is the analytic continuation of $\arctan{x}$ from the real line to $\arctanDom$.
Similarly, the principal branch $\arccot{z}$ is the analytic continuation of $\arccot{x}$ from $\mathbb{R}\setminus\{0\}$ to $\arccotDom$.
On the $y$-axis away from the cuts, taking the limit as $x\to 0$ shows that
$$%\begin{equation}
\begin{aligned}
\arctan(\rmi y) &= \rmi\arctanh{y}  & (y^2<1), \\
\arccot(\rmi y) &= -\rmi\arccoth{y} & (y^2>1).
\end{aligned}
$$%\end{equation}
On the cuts, $\arctan{z}$ and $\arccot{z}$ are two-valued,
$$%\begin{equation}
\begin{aligned}
\arctan(\pm 0+\rmi y) &= \pm\frac{\pi}{2}+\rmi\arccoth{y} & (y^2>1), \\
\arccot(\pm 0+\rmi y) &= \pm\frac{\pi}{2}-\rmi\arctanh{y} & (y^2<1).
\end{aligned}
$$%\end{equation}
The first equation here is expressed in terms of the natural logarithm in the DLMF\cite[Eq.~4.23.27]{DLMF}.
We see that the cut planes $\arctanDom$ under $\arctan{z}$ and $\arccotDom$ under $\arccot{z}$ both map to $[-\pi/2,\pi/2]\times\mathbb{R}$.
We verify that
$$%\begin{equation}
\begin{aligned}
\frac{\rmd}{\rmd z}\arctan{z} &= \frac{1}{1+z^2},  & \int\arctan{z}\,\rmd z &= z\arctan{z}-\frac{1}{2}\ln(1+z^2)+C, \\
\frac{\rmd}{\rmd z}\arccot{z} &= -\frac{1}{1+z^2}, & \int\arccot{z}\,\rmd z &= z\arccot{z}+\frac{1}{2}\ln(1+z^2)+C.
\end{aligned}
$$%\end{equation}

\section{Principal branch of the inverse hyperbolic sine}
The inverse hyperbolic sine function has branch points at $\pm\rmi$ and $\infty$\,\cite{Wolfram}.
The princpal branch $\arcsinh{z}$ is defined on the cut plane $\arcsinhDom$ by\cite[Eq.~4.37.1]{DLMF}
$$%\begin{equation}
\arcsinh{z}=\int_{0}^{z}\frac{\rmd t}{\sqrt{1+t^2}}=-\rmi\int_{0}^{\rmi z}\frac{\rmd t}{\sqrt{1-t^2}}=\rmi\arcsin(-\rmi z).
$$%\end{equation}
From the equation above for $\arcsin{z}$ we discover that
$$%\begin{equation}
\arcsinh{z} = \sign(x)\arccosh{a}+\rmi\arcsin{b},
$$%\end{equation}
where $a=(\rho_1+\rho_2)/2$ and $b=(\rho_2-\rho_1)/2$ with $a\ge1$, $b\in[-1,1]$ ($b>0$ for $y>0$ and $b<0$ for $y<0$)
Here $\rho_1=|z-\rmi|=\sqrt{x^2+(y-1)^2}$ and $\rho_2=|z+\rmi|=\sqrt{x^2+(y+1)^2}$.
Since $\rho_1(-z)=\rho_2(z)$ and $\rho_2(-z)=\rho_1(z)$,  $a(-z)=a(z)$ and $b(-z)=-b(z)$.
We note that the partial derivatives of $a$ and $b$ with respect to $y$ are
$$%\begin{equation}
a_y = \frac{b(a^2-1)}{a^2-b^2}, \quad b_y = \frac{a(1-b^2)}{a^2-b^2}.
$$%\end{equation}
We observe that $\arcsinh(-z)=-\arcsinh{z}$.
This reflection formula is in the DLMF\cite[Eq.~4.37.10]{DLMF}.
On the $x$-axis, $a=\sqrt{1+x^2}$, $b=0$ and we confirm that $\arcsinh{z} = \arcsinh{x}$.
The principal branch $\arcsinh{z}$ is the analytic continuation of $\arcsinh{x}$ from the real line to the cut plane $\arcsinhDom$.
On the $y$-axis away from the cut, $z=\rmi y\in[-\rmi,\rmi]$, $a=1$, $b=y$ and
$$%\begin{equation}
\arcsinh(\rmi y)=\rmi\arcsin{y} \quad (y^2\le1).
$$%\end{equation}
On the $y$-axis on the cut, $z=\pm0+\rmi y (y^2>1)$, $a=|y|$, $b=\sign(y)$, and $\arcsinh{z}$ is two-valued,
$$%\begin{equation}
\arcsinh(\pm0+\rmi y) = \pm\arccosh|y|+\rmi\sign(y)\frac{\pi}{2} \quad (y^2>1).
$$%\end{equation}
The corresponding equations in the DLMF are\cite[Eq.~4.37.17--18]{DLMF}.
We see that $\arcsinh{z}$ maps the cut plane $\arcsinhDom$ to $\mathbb{R}\times[-\pi/2,\pi/2]$.

Working backwards, we infer that\cite[Eq.~4.37.1]{DLMF}
$$%\begin{equation}
\arcsinh{z} = \int_{0}^{z}\frac{\rmd t}{\sqrt{1+t^2}} = \int_0^z\frac{\rmd t}{\Z(t)},
$$%\end{equation}
where
$$%\begin{equation}
\Z(z)=a\sqrt{1-b^2}+\rmi\sign(x)b\sqrt{a^2-1},
$$%\end{equation}
is the analytic continuation of $\sqrt{1+x^2}$ from the real line to the cut plane $\arcsinhDom$.
We note that $|\Z(z)|=\sqrt{\rho_1\rho_2}$ and $\Z(-z)=\Z(z)$.
On the $x$-axis, $a=\sqrt{1+x^2}$, $b=0$ and we confirm that
$$%\begin{equation}
\Z(x) = \sqrt{1+x^2} \quad (x\in\mathbb{R}).
$$%\end{equation}
On the $y$-axis away from the cut, $a=1$, $b=y$ and
$$%\begin{equation}
\Z(\rmi y) = \sqrt{1-y^2} \quad (y^2\le 1).
$$%\end{equation}
On the $y$-axis on the cut, $a=|y|$, $b=\sign(y)$ and $\Z(z)$ is two-valued,
$$%\begin{equation}
\Z(\pm0+\rmi y) =+0 \pm\rmi\sign(y)\sqrt{y^2-1} \quad (y^2>1).
$$%\end{equation}
We see that $\Z(z)$ maps the cut plane $\arcsinhDom$ to the closed right half plane.
Away from the axes, the first and third quadrants map under $\Z(z)$ to the first quadrant while the second and fourth quadrants map to the fourth quadrant.
We verify that $\Z(z)^2=1+z^2$, $\Z'(z)=z/\Z(z)$, and
$$%\begin{equation}
\frac{\rmd}{\rmd z}\arcsinh{z} = \frac{1}{\Z(z)}, \quad \int\arcsinh{z}\,\rmd z = z\arcsinh{z}-\Z(z)+C.
$$%\end{equation}

\section{Principal branch of the inverse hyperbolic cosine}
The inverse hyperbolic cosine function has branch points at $\pm1$ and $\infty$\,\cite{Wolfram}.
The principal branch $\arccosh{z}$ is defined on the cut plane $\arccoshDom$ by\cite[Eq.~4.37.2]{DLMF}
$$%\begin{equation}
\arccosh{z} = \int_{1}^{z}\frac{\rmd t}{\sqrt{t^2-1}} = \int_1^z\frac{\rmd t}{\X(t)}.
$$%\end{equation}
If we let $c(z)=\arccosh{\alpha}+\rmi\sign(y)\arccos{\beta}$ we find that $c'(z)=1/\X(z)$. It follows that
$$%\begin{equation}
\arccosh{z} = c(z)-c(1) = \arccosh{\alpha}+\rmi\sign(y)\arccos{\beta}.
$$%\end{equation}
This expression is in Johnson\cite[p.~264]{Johnson}.
For $z$ in the cut plane, if $z\notin[1,\infty)$ then $-z$ is also in the cut plane and $\arccosh(-z)=-\rmi\sign(y)\pi+\arccosh{z}$.
This reflection formula has $-\rmi\pi$ for $y>0$ and $+\rmi\pi$ for $y<0$.
The current release of the DLMF\cite[Eq.~4.37.11]{DLMF} erroneously has $+\rmi\pi$ for $y>0$ and $-\rmi\pi$ for $y<0$.
On the $x$-axis away from the cut, $x\ge1$, $\alpha=x$, $\beta=1$ and we confirm that $\arccosh{z}=\arccosh{x}$.
The principal branch $\arccosh{z}$ is the analytic continuation of $\arccosh{x}$ from $[1,\infty)$ to the cut plane $\arccoshDom$.
On the $x$-axis on the cut (for $x^2<1$, $\alpha=1$ and $\beta=x$ while for $x<-1$, $\alpha=-x$ and $\beta=-1$), $\arccosh{z}$ is two-valued,
$$%\begin{equation}
\arccosh(x\pm\rmi 0)
=
\begin{cases}
\hfill\pm\rmi\arccos{x}, & \hfill (x^2<1), \\
\arccosh(-x)\pm\rmi\pi, & (x<-1).
\end{cases}
$$%\end{equation}
The corresponding equations in the DLMF are \cite[Eq.~4.37-22--23]{DLMF}.
On the $y$-axis, $z=\rmi y$, $\alpha=\sqrt{1+y^2}$, $\beta=0$, and
$$%\begin{equation}
\arccosh(\rmi y) = \arcsinh{|y|}+\rmi\sign(y)\frac{\pi}{2} \quad (y^2>0).
$$%\end{equation}
The corresponding equation in the DLMF is \cite[Eq.~4.37.20]{DLMF}.
We see that $\arccosh{z}$ maps the cut plane $\arccoshDom$ to $[0,\infty)\times[-\pi,\pi]$.
We verify that
$$%\begin{equation}
\frac{\rmd}{\rmd z}\arccosh{z} = \frac{1}{\X(z)}, \quad  \int\arccosh{z}\,\rmd z = z\arccosh{z}-\X(z)+C.
$$%\end{equation}

\section{Principal branches of the inverse hyperbolic cosecant and inverse hyperbolic secant}
The inverse hyperbolic cosecant function has branch points at $\pm\rmi$ and $0$\,\cite{Wolfram}.
The principal branch $\arccsch{z}$ is defined on the cut plane $\arccschDom$ by\cite[Eq.~4.37.7]{DLMF}
$$%\begin{equation}
\begin{aligned}
\arccsch{z} &= \arcsinh(1/z).
\end{aligned}
$$%\end{equation}
It follows that $\arccsch(-z)=-\arccsch{z}$.
This reflection formula is in the DLMF\cite[Eq.~4.37.13]{DLMF}.
On the $x$-axis away from the cut, $x^2>0$ and $\arccsch{z}=\arccsch{x}$.
The principal branch $\arccsch{z}$ is the analytic continuation of $\arccsch{x}$ from $\mathbb{R}\setminus\{0\}$ to $\arccschDom$.
On the $y$-axis away from the cut,
$$%\begin{equation}
\arccsch(\rmi y) = -\rmi\arccsc{y} \quad (y^2\ge1).
$$%\end{equation}
On the $y$-axis on the cut, $\arccsch{z}$ is two-valued,
$$%\begin{equation}
\arccsch(\pm0+\rmi y) = \pm\arcsech|y|-\rmi\sign(y)\frac{\pi}{2} \quad (0<y^2<1).
$$%\end{equation}
We see that $\arccsch{z}$ maps the cut plane $\arccschDom$ to $\mathbb{R}\times[-\pi/2,\pi/2]$.

The inverse hyperbolic secant function has branch points at $\pm1$ and $0$\,\cite{Wolfram}.
The principal branch $\arcsech{z}$ is defined on the cut plane $\arcsechDom$ by\cite[Eq.~4.37.8]{DLMF}
$$%\begin{equation}
\arcsech{z} = \arccosh(1/z).
$$%\end{equation}
For $z$ in the cut plane, if $z\notin(0,1]$ then $-z$ is also in the cut plane and $\arcsech(-z)=\rmi\sign(y)\pi+\arcsech{z}$.
This reflection formula has $+\rmi\pi$ for $y>0$ and $-\rmi\pi$ for $y<0$.
The current release of the DLMF\cite[Eq.~4.37.14]{DLMF} erroneously has $-\rmi\pi$ for $y>0$ and $+\rmi\pi$ for $y<0$.
On the $x$-axis away from the cut, $x\in(0,1]$ and $\arcsech{z}=\arcsech{x}$.
The principal branch $\arcsech{z}$ is the analytic continuation of $\arcsech{x}$ from $(0,1]$ to the cut plane $\arcsechDom$.
On the $x$-axis on the cut, $\arcsech{z}$ is two-valued,
$$%\begin{equation}
\arcsech(x\pm\rmi 0)
=
\begin{cases}
\hfill\mp\rmi\arcsec{x}, & \hfill (x^2>1), \\
\arcsech(-x)\mp\rmi\pi, & (-1<x<0).
\end{cases}
$$%\end{equation}
On the $y$-axis,
$$%\begin{equation}
\arcsech(\rmi y) = \arccsch|y|-\rmi\sign(y)\frac{\pi}{2} \quad (y^2>0).
$$%\end{equation}
We see that $\arcsech{z}$ maps the cut plane $\arcsechDom$ to $\mathbb{R}\times[-\pi,\pi]$.
We verify that
$$%\begin{equation}
\begin{aligned}
\frac{\rmd}{\rmd z}\arccsch{z} &= -\frac{1}{z^2\Z(1/z)}, & \int\arccsch{z}\,\rmd z &= z\arccsch{z}+\arctanh(\Z(1/z))+C, \\
\frac{\rmd}{\rmd z}\arcsech{z} &= -\frac{1}{z^2\X(1/z)}, & \int\arcsech{z}\,\rmd z &= z\arcsech{z}-\arctan(\X(1/z))+C.
\end{aligned}
$$%\end{equation}

\section{Principal branches of the inverse hyperbolic tangent and inverse hyperbolic cotangent}
The inverse hyperbolic tangent and inverse hyperbolic cotangent functions have branch points at $\pm1$\,\cite{Wolfram}.
The principal branch $\arctanh{z}$ of the inverse hyperbolic tangent is defined on the cut plane $\arctanhDom$ by\cite[Eq.4.37.3]{DLMF}
$$%\begin{equation}
\arctanh{z}=\int_0^z\frac{\rmd t}{1-t^2}=-\rmi\int_0^{\rmi z}\frac{\rmd t}{1+t^2}=-\rmi\arctan(\rmi z).
$$%\end{equation}
On the $x$-axis away from the cut, $x^2<1$ and we confirm that $\arctanh{z}=\arctanh{x}$.
The principal branch $\arctanh{z}$ is the analytic continuation of $\arctanh{x}$ from $(-1,1)$ to the cut plane $\arctanhDom$.

The principal branch $\arccoth{z}$ of the inverse hyperbolic cotangent is defined on the cut plane $\arccothDom$ by\cite[Eq.~4.37.9]{DLMF}
$$%\begin{equation}
\arccoth{z}=\arctanh(1/z)=-\rmi\arctan(\rmi/z)=\rmi\arccot(\rmi z).
$$%\end{equation}
On the $x$-axis away from the cut, $x^2>1$ and we confirm that $\arccoth{z}=\arccoth{x}$.
The principal branch $\arccoth{z}$ is the analytic continuation of $\arccoth{x}$ from $(-\infty,-1)\cup(1,\infty)$ to the cut plane $\arccothDom$.

It follows that $\arctanh(-z)=-\arctanh{z}$ and $\arccoth(-z)=-\arccoth{z}$.
These reflection formulae are in the DLMF\cite[Eq.~4.37.12,15]{DLMF}.
The reflection formula for $\arctanh{z}$,
$$%\begin{equation}
\arctanh{z}=\rmi\sign(y)\frac{\pi}{2}+\arccoth{z} \quad (y\ne 0),
$$%\end{equation}
follows from the reflection formula for $\arctan{z}$.
From the expression for $\arccot{z}$,
$$%\begin{equation}
\arccoth{z}=\frac{1}{4}\ln\frac{(x+1)^2+y^2}{(x-1)^2+y^2}+\frac{\rmi}{2}\left(\arctan\frac{x-1}{y}-\arctan\frac{x+1}{y}\right) \quad (y\ne 0).
$$%\end{equation}
On the cuts, $\arctanh{z}$ and $\arccoth{z}$ are two-valued,
$$%\begin{equation}
\begin{aligned}
\arctanh(x\pm\rmi 0) &= -\rmi\arctan(\mp 0+\rmi x) = \arccoth{x}\pm\rmi\frac{\pi}{2} & (x^2>1), \\
\arccoth(x\pm\rmi 0) &= \rmi\arccot(\mp 0+\rmi x)  = \arctanh{x}\mp\rmi\frac{\pi}{2} & (x^2<1).
\end{aligned}
$$%\end{equation}
On the imaginary axis,
$$%\begin{equation}
\begin{aligned}
\arctanh(\rmi y) &= \rmi\arctan{y}  & (y\in\mathbb{R}), \\
\arccoth(\rmi y) &= -\rmi\arccot{y} & (y^2>0).
\end{aligned}
$$%\end{equation}
We see that the cut planes $\arctanhDom$ under $\arctanh{z}$ and $\arccothDom$ under $\arccoth{z}$ both map to $\mathbb{R}\times[-\pi/2,\pi/2]$.
We verify that
$$%\begin{equation}
\begin{aligned}
\frac{\rmd}{\rmd z}\arctanh{z} &= \frac{1}{1-z^2}, & \int\arctanh{z}\,\rmd z &= z\arctanh{z}+\frac{1}{2}\ln(z^2-1)+C, \\
\frac{\rmd}{\rmd z}\arccoth{z} &= \frac{1}{1-z^2}, & \int\arccoth{z}\,\rmd z &= z\arccoth{z}+\frac{1}{2}\ln(z^2-1)+C.
\end{aligned}
$$%\end{equation}

\section{Concluding Remarks}
The principal branch expressions for $\sqrt{z^2-1}$, $\sqrt{1-z^2}$ and $\sqrt{1+z^2}$ (here denoted by $\X(z)$, $\Y(z)$ and $\Z(z)$) are new 
breakthrough discoveries that lead smoothly to principal branch expressions for $\arcsin{z}$, $\arccos{z}$, $\arcsinh{z}$ and $\arccosh{z}$.
The principal branch expressions for $\arcsin{z}$, $\arccos{z}$ and $\arccosh{z}$ have long been in print in Appendix B of Johnson\cite{Johnson}.
The latter, first published in 1924, was based on a popular set of in-house course notes on electrical communication at Bell Telephone Laboratories.
Abramowitz and Stegun\cite{A&S}, published four decades later and a favorite of the scientific community, is not complete.
Only the general principal branch expressions for $\arcsin{z}$ and $\arccos{z}$, and the reflection formula binding $\arctan{z}$ and $\arccot{z}$, are in the DLMF to this day.
Using this article, the general principal branch expressions for $\arccot{z}$, $\arcsinh{z}$, $\arccosh{z}$, and $\arccoth{z}$,
and the reflection formula binding $\arctanh{z}$ and $\arccoth{z}$, could be added to the DLMF.
The axis traces, derivatives, and antiderivatives, could also be added.

\end{document}